\documentclass[notitlepage,leqno,11pt]{article}
\usepackage{amssymb,cite,amsthm}
\catcode`\@=11 \@addtoreset{equation}{section}

\catcode`\@=12

\usepackage{amsrefs}
\usepackage{latexsym}
\usepackage{textcomp}
\usepackage{amsmath}
\usepackage{amsfonts}
\usepackage{amssymb}
\usepackage{mathrsfs}
\usepackage{dsfont}

\renewcommand{\a }{\alpha }
\renewcommand{\b }{\beta }
\renewcommand{\d}{\delta }

\newcommand{\D }{\Delta }

\newcommand{\e }{\varepsilon }
\newcommand{\g }{\gamma}

\renewcommand{\l }{\lambda }
\renewcommand{\L }{\Lambda }

\newcommand{\n }{\nabla }
\newcommand{\vp }{\varphi }
\renewcommand{\phi}{\varphi}

\newcommand{\s }{\sigma }

\renewcommand{\th }{\theta }

\renewcommand{\o }{\omega }

\newcommand{\be}{\begin{equation}}
\newcommand{\ee}{\end{equation}}

\newcommand{\R}{\mathbb{R}}
\newcommand{\N}{\mathbb{N}}

\newcommand{\de}{\partial}

\newcommand{\ra}{{\rangle}}
\newcommand{\la}{{\langle}}

\renewcommand{\k}{\kappa}

\newcommand{\calH }{\mathcal{H}}

\newcommand{\calB }{\mathcal{B}}

\newcommand{\calN}{{\mathcal N}}
\newcommand{\calJ}{\mathcal{ J}}
\newcommand{\calM}{{\mathcal M}}

\renewcommand{\le}{\leqslant}
\renewcommand{\leq}{\leqslant}
\renewcommand{\ge}{\geqslant}
\renewcommand{\geq}{\geqslant}

\newtheorem{Theorem}{Theorem}[section]
\newtheorem{Lemma}[Theorem]{Lemma}
\newtheorem{Proposition}[Theorem]{Proposition}
\newtheorem{Corollary}[Theorem]{Corollary}

\def\proof{\noindent{{\bf Proof. }}}
\def\square{\vbox{
    \hrule height .4pt
    \hbox{\vrule width .4pt height 7pt \kern 7pt
       \vrule width .4pt}
    \hrule height .4pt }}

\def\square{\vbox{
    \hrule height .4pt
    \hbox{\vrule width .4pt height 7pt \kern 7pt
       \vrule width .4pt}
    \hrule height .4pt }}

\def\QED{\hfill {$\square$}\goodbreak \medskip}

\def\R{{\mathbb R}}

\def\div{{\rm div}}

\newcommand{\Ds}{(-\Delta)^s}

\newcommand{\Dsib}{(-\Delta)^{\bar{\s}}}

\def\RNp{{\mathbb R}^{N+1}_+}
\font\sc=cmcsc9  \textwidth=14truecm
\title{Uniqueness and nondegeneracy of positive solutions \\
of $\Ds u+u=u^p$ in $\R^N$\\
when $s$ is close to $1$}
\author{Mouhamed Moustapha Fall and Enrico Valdinoci}

\begin{document}
\date{}
\maketitle

 \let\thefootnote\relax\footnotetext{
M. M.~Fall,
African Institute for Mathematical Sciences (A.I.M.S.) of Senegal,
KM 2, Route de Joal (Centre I.R.D. Mbour),
B.P. 1418 Mbour, S\'en\'egal,
{\tt{mouhamed.m.fall@aims-senegal.org}}.
E.~Valdinoci,
Weierstra{\ss} Institut f\"ur Angewandte Analysis und Stochastik (W.I.A.S.),
Mohrenstra{\ss}e 39,
10117 Berlin, Germany,
{\tt{enrico.valdinoci@wias-berlin.de}}.}

\bigskip

\begin{abstract}
We consider the equation $\Ds u+u=u^p$, with $s\in(0,1)$
in the subcritical range of $p$.
We prove that if $s$ is sufficiently close to $1$
the equation possesses a unique minimizer, which is nondegenerate.
\end{abstract}

\section{Introduction}

The purpose of this paper is to provide some
nondegeneracy and uniqueness result for solutions
of an equation driven by a nonlocal operator.
In striking contrast with the local case, extremely
little is known about these topics in the nonlocal framework
and a satisfactory analysis of the problem is still largely missing,
in spite of some striking recent contributions in specific cases.

Our approach is to obtain
some nondegeneracy and uniqueness results by compactness
and bifurcation arguments
from the local case, that is when the fractional
parameter involved is sufficiently close to being an integer.

Let us introduce the setting in which the problem is posed.
Let $N\ge2$ be the dimension of the ambient space $\R^N$
and let $s\in(0,1]$
be our fractional parameter.

We consider the fractional exponent
$$ 2^*_s:=\left\{
\begin{matrix}
\displaystyle\frac{2N}{N-2s} & {\mbox{ if $N\ge3$, or $N=2$ and $s\in(0,1)$,}}\\
+\infty & {\mbox{ if $N=2$ and $s=1$.}}
\end{matrix}\right.$$
We recall that
this exponent plays the role of the classical critical Sobolev exponent
for the fractional Sobolev spaces (see, e.g., \cite{guida} for a gentle
introduction to the topic, and notice that~$2^*_s$ is increasing in~$s$
and coincides with the classical Sobolev exponent for~$s=1$). We consider here the fractional
Sobolev space
$$
H^s(\R^N):=\Big\{ u\in L^2(\R^N)\,:\,\int_{\R^N}|\xi|^{2s} |\widehat{u}|^2\,d\xi<\infty  \Big\},
$$
with norm
$$
\|u\|_s^2:=\int_{\R^N}(1+|\xi|^{2s}) |\widehat{u}|^2\,d\xi=
\|u\|^2_{L^{2}(\R^N)}+\int_{\R^N} |\xi|^{2s} |\widehat{u}|^2\,d\xi,
$$
where, as usual, $\widehat{u}$ is the Fourier transform of the function $u$,
namely
$$
\widehat{u}(\xi):=\frac{1}{(2\pi)^{\frac{N}{2}}}\,
\int_{\R^N}e^{-\imath\xi\cdot{x} }u(x)dx.
$$
We also denote by $H^s_{rad}(\R^N)$ the space of the radially
symmetric functions of $H^s(\R^N)$.
We recall that $2^*_s$ provides a compactness threshold
for such radial functions, since
$L^{q}(\R^N) $
is compactly embedded in $H^s_{rad}(\R^N)$ for every $q\in (1,2^*_s)$
(see Proposition 1.1 in~\cite{lions}).

In this functional framework, we are concerned
with the uniqueness and nondegeneracy properties of the positive functions
solving the fractional elliptic semilinear problem
\be
\label{eq:pblm}
\Ds u+u =u^p \quad\textrm{ in }\R^N.
\ee
Here we take $p\in (1,2^*_s-1)$ (i.e., the exponent
$p+1$ is subcritical with respect to the above mentioned embeddings).
Problems of this type has received a great attention recently,
both by themselves and in connection with solitary solutions
of nonlinear dispersive wave equations
(such as the Benjamin-Ono equation, the Benjamin-Bona-Mahony equation
and the fractional Schr\"odinger equation, see e.g. \cites{W87, bona97, ann97, maris, martel, FL}).

In this framework, the classical, local Hamiltonian operator is replaced
by a fractional, nonlocal one, and the classical diffusion induced by
Brownian motions is replaced by a non-local diffusion driven by
$2s$-stable L\'evy processes.

These type of fractional operators are now becoming also very popular in real-world models
(for instance in financial mathematics, nonlocal
stochastic control,
nonlocal electrostatics, denoising and image processing,
oceanography, 
dislocation dynamics in crystals, etc.), see for instance~\cite{guida}
and references therein.

Since the fractional Laplacian of $\vp\in C^\infty_c(\R^N)$ may be defined via Fourier transform as
\be\label{FOU}
 \widehat {\Ds \vp}(\xi):=|\xi|^{2s} \widehat \vp (\xi) \qquad \text{for
  $\xi \in \R^N$,}
\ee
we may apply
Plancherel's formula and adopt a weak (or distributional)
notion of solution $u\in H^s(\R^N)$ for problem \eqref{eq:pblm} via the identity
$$
\frac12
\int_{\R^N} |\xi|^{2s}(\widehat{u}\,\overline{\widehat{\vp}}+
\overline{\widehat{u}}\,{\widehat{\vp}}) \,d\xi=
\int_{\R^N} |\xi|^{2s}\widehat{u}\,\overline{\widehat{\vp}} \,d\xi=\int_{\R^N}(u^p-u)\vp\,dx$$
for any $\vp\in H^s(\R^N)$.
This notion of solution may be reduced to the one in the viscosity sense
(see \cites{servadei-reg, serra-reg}) and therefore the fractional
Laplace regularity theory applies (see \cite{Sil}).
It is known that problem \eqref{eq:pblm} admits a positive radial
solution (see \cites{dipierro, FQT}). Such solution is called a ground state,
since it is obtained (up to scaling) by a constrained minimization
problem of the functional
$$
J_s(u,\nu):=\frac{1}{2}\|u\|_s^2-\frac{\nu}{p+1} \int_{\R^N}|u|^{p+1}\,dx,
$$
namely it attains the following greatest lower bound:
\be\label{eq:nus}
\nu_s:=\inf_{u\in H^s(\R^N)}\frac{\|u\|^2_s}{
\displaystyle\left(\int_{\R^N}|u|^{p+1}\right)^{2/(p+1)}}
= \inf_{\genfrac{}{}{0pt}{}{\scriptstyle{u\in H^s(\R^N)}}{\scriptstyle{ \|u\|_{L^{p+1}(\R^N)}=1 }}} \|u\|^2_s.
\ee
We observe that if $u_s$ is such that $\|u_s\|_{L^{p+1}(\R^N)}=1$
and $\nu_s= \|u\|^2_s$, than it is a solution of
\begin{equation}\label{EQU} \Ds u_s+u_s =\nu_s u_s^p\end{equation}
and so it solves \eqref{eq:pblm} (up to scaling). Also its derivatives
$\partial_i u_s$ are solution of the linearized equation
\be\label{LIN}
\Ds (\partial_i u_s)+\partial_i u_s =p\nu_s u_s^{p-1} \partial_i u_s\ee
and therefore \be\label{1.4bis}
{\mbox{$\partial_i u_s$ belongs to the kernel of the operator $J_s''(u_s,\nu_s)$.}}\ee

The first result of this paper is nondegeneracy, namely
that these derivative and their linear combinations exhaust
$Ker(J_s''(u_s,\nu_s))$ at least when
$s$ is sufficiently close to $1$
(of course, since
we are interested here in the case $s$ close to $1$
with a fixed exponent $p$,
we fix~$S\in(0,1)$ and $p \in 2^*_{S} - 1$,
and all the arguments we present assume implicitly
that~$s\in [S,1]$).

\begin{Theorem}\label{th:nondegen}
There exists $s_0\in(0,1)$ such that for every $s\in(s_0,1)$
if $u_s$ is a minimizer for $\nu_s$ then
$$
Ker(J_s''(u_s,\nu_s))=span\{\de_i u_s,\,i=1,\dots,N\}.
$$
\end{Theorem}

Our next result is a uniqueness property.

\begin{Theorem}\label{th:Uniq}
There exists $s_0\in(0,1)$ such that for every $s\in(s_0,1)$, the minimizer for $\nu_s$ is unique, up to translations.
\end{Theorem}


In the local case $s=1$, the results in Theorems \ref{th:nondegen}
and \ref{th:Uniq} were obtained in \cites{serrin-mcleod, mcleod, kwong}
but the specific arguments used there
are not directly applicable to the nonlocal case $s\in(0,1)$.
Before this paper, the only results
available in the nonlocal case were the ones obtained in
\cite{AT} for $N=1$, $s=1/2$ and $p=2$, and recently extended in
\cite{FL} for $N=1$ and all $s\in(0,1)$.

After this paper was completed, arxived in \cite{FV} and submitted,
the striking paper \cite{FLS} has appeared, showing that
Theorems \ref{th:nondegen}
and \ref{th:Uniq} hold for any~$s\in(0,1)$.

We also point out that, soon after \cite{FV},
some interesting nondegeneracy results have been obtained in \cite{DPS}
for a related, but
different, fractional problem.

For other recent variational problems related to the fractional
Laplacian see, for instance, \cites{servadei, servadeiDCDS, FW, secchi} and
references therein.
The rest of the paper is organized as follows. In Section \ref{S:1}
we collect some preliminary material, likely well-known to the expert
readers, concerning some uniform estimates on the minimizers,
some related asymptotics and a (up to now classical)
local realization of the fractional Laplacian. Then, in
Section~\ref{S:2}, we prove the nondegeneracy result of
Theorem~\ref{th:nondegen}. The uniqueness result of
Theorem~\ref{th:Uniq} is proved in Sections~\ref{S:3} and~\ref{S:4},
by combining a series of arguments related  to the construction of a branch of
pseudo-minimizers $U_1+\o_s$, with $s$ varies near 1, which are uniquely determined by their perturbation $\o_s$.
Uniqueness is then deduced by showing that
radially symmetric minimizers belongs to such a branch.


\section{Preliminaries}\label{S:1}

\subsection{Uniform estimates and asymptotics}
By Lion's concentration compactness, minimizers for $\nu_s$ always exists
(see, e.g.~\cites{dipierro, FQT} for details) and do not change sign. In this paper, we will consider only positive minimizers. They are radially symmetric by \cite{FQT}
(and, as usual, we take the center of symmetry to be the origin of~$\R^N$).
The minimizers attain the minimal value~$\nu_s$ of the functional in~\eqref{eq:nus}
and they are normalized to have norm~$1$
in~$L^p(\R^N)$. Also, thanks to Theorem 1.2 in \cite{FQT}, we have the decay estimate
$$
u_s\leq C |x|^{-(N+2s)}\quad\textrm{ in } \R^N.
$$
We call $\calM_{s}$ the the space of these positive, radially symmetric even minimizers $u_s$ for $\nu_s$
normalized so that $\|u_s\|_{L^{p+1}(\R^N)}=1$.
Therefore if $u_s\in \calM_{s}$
then
\begin{equation}\label{radially}
\|u_s\|_{L^\infty(\R^N)}=|u_s(x_0^s)|,
\end{equation}
for some $x_0^s\in \R^N$.
Now we state a uniform bound on~$\nu_s$:

\begin{Lemma}\label{88}
We have that~$\displaystyle\sup_{s\in(0,1]}\nu_s<+\infty$.
\end{Lemma}

\proof Let~$u_1\in\calM_1$.
Notice that~$ |\xi|^{2s}\le 1+|\xi|^2$ and therefore
$$ \|u_1\|_s \le
2\|u_1\|_{L^{2}(\R^N)}+\int_{\R^N} |\xi|^{2} |\widehat{u}|^2\,d\xi\le 2\|u_1\|_1=2\nu_1.$$
Since~$\nu_s\le \|u_1\|_s$, the desired result follows.
\QED

The following result provides uniform bounds on the minimizers.

\begin{Lemma}\label{lem:unif}
Given $s_0\in(0,1)$, we have
\be\label{eq:unif}
0<\g_{s_0}:=\sup_{s\in(s_0,1)}\sup_{u_s\in \calM_{s} }\|u_s\|_{L^\infty(\R^N)}<\infty.
\ee
Also, given $s_1>1/2$ and $\b\in(0,1)$,
\be\label{eq:unif-grad}
\sup_{s\in(s_1,1)}\sup_{u_s\in \calM_{s}}\|u_s\|_{C^{1,\b}(\R^N)}<\infty.
\ee
\end{Lemma}
\proof The first inequality in \eqref{eq:unif} is obvious since
$$ \g_{s_0}\ge \sup_{u_1\in \calM_{1} }\|u_1\|_{L^\infty(\R^N)}>0.$$
Now we prove the second inequality in \eqref{eq:unif}. For this, we
define
\begin{equation}\label{7878}
\l_s:=\|u_s\|_{L^\infty(\R^N)}
\end{equation}
and we argue by contradiction: we
suppose that $\l_s\to \infty$ for a sequence $s\to\bar\s\in[s_0,1]$.
We set
$$
v_s(x):=\l_s^{-1}u_s(\l_s^{\frac{2}{2s-N}}x+x_0^s)
$$
so that
$$
\|v_s\|_{L^\infty(\R^N)}=1=v_s(0),
$$
$$\widehat{v_s}(\xi)=\lambda_s^{-1+\frac{2N}{N-2s}}e^{\imath\xi\cdot x_0^s}\widehat{u_s}(\lambda_s^{\frac{2}{N-2s}}\xi)$$
and,
by Lemma \ref{88},
$$
\int_{\R^N}|\xi|^{2s}|\widehat{v_s}|^2\,d\xi=
\int_{\R^N}|\xi|^{2s}|\widehat{u_s}|^2\,d\xi\leq
\nu_s\leq {\mbox{Const}}.$$
Therefore  $v_s \rightharpoonup v$ in $H^t(\R^N)$ for every $t<\bar{\s}$ and
$$
{\mbox{$v_s\to v$ in $L^{2}_{loc}(\R^N)$}}.
$$
Also, from \eqref{EQU},
\begin{equation}\label{HH:0}
\Ds v_s(x)=-\l_s^{\frac{4s}{-N+2s}}v_s(x)+ \l_s^{p-\frac{N+2s}{N-2s}}\nu_s
v_s^p(x).
\end{equation}
Now we recall Proposition 2.1.9 in~\cite{Sil},
according to which
we have that there is a constant $C(s,N,\a)$  such that
\begin{equation}\label{HH:1}
\|v_s\|_{C^{0,\a}(\R^N)}\leq C(s,N,\a)\left( \|\Ds v_s\|_{L^\infty(\R^N)} + \| v_s\|_{L^\infty(\R^N)}  \right),
\end{equation}
where one can fix~$\a< 2\bar\s$ for $2\bar\s<1$ and $\a< 2\bar\s-1$ for $2\bar\s>1$
and the constant $C(s,N,\a)$ is bounded uniformly in $s\in[s_0,1]$. {F}rom
Lemma \ref{88},
\eqref{HH:0} and
\eqref{HH:1}, we see that
$\|v_s\|_{C^{0,\a}(\R^N)}$ is bounded uniformly when $s\to\bar\s$.
Accordingly, by the Ascoli theorem, we may suppose that $v_s$ converges
locally uniformly to $v$ and passing to the limit in \eqref{HH:0},   we have that $v\equiv 0$. In particular
$$ 0=\lim_{s\to\bar\s} |v_s(0)|=\lim_{s\to\bar\s} \lambda_s^{-1} |u_s(x_0^s)|=1,$$
due to \eqref{radially}
and \eqref{7878}.
This is a contradiction and so \eqref{eq:unif}
is proved.

To prove \eqref{eq:unif-grad}
we use once again
Proposition 2.1.9 in~\cite{Sil}, see also \cite{CS}, according to which, for any $s\in(s_1,1]$,
$$
\|u_s\|_{C^{1,\b}(\R^N)}\leq C(s,N,\a)\left( \|\Ds u_s\|_{L^\infty(\R^N)} + \| u_s\|_{L^\infty(\R^N)}  \right),
$$
where $C(s,N,\a) $ is uniformly bounded on $[s_1,1]$.
Then, the latter inequality implies \eqref{eq:unif-grad},
thanks to \eqref{EQU},
\eqref{eq:unif} and
Lemma \ref{88}.
\QED

\begin{Corollary}\label{C:REG}
Given $s_0\in(0,1)$, we have
$$
\sup_{s\in(s_0,1)}\sup_{u_s\in \calM_{s} }
\|u_s\|_{2s}
<\infty.$$
\end{Corollary}

\proof Let $s_0\in(0,1)$,
$u_s\in\calM_s$ and $f_s(x):=\nu_s
u_s^p(x)-u_s(x)$.
Notice that
$$ \int_{\R^N}|u_s|^{2p}\,dx\le \|u_s\|_{L^\infty(\R^N)}^{2(p-1)}
\int_{\R^N}|u_s|^{2}\,dx\le C_1,$$
with $C_1>0$ independent of $s$ and $u_s$, thanks to
\eqref{eq:unif},
Lemma \ref{88} and the fact that $p>1$. Moreover,
$$ \|u_s\|_{L^2(\R^N)}^2\le \nu_s\le C_2,$$
with $C_2>0$ independent of $s$ and $u_s$, thanks to
Lemma \ref{88}. As a consequence, and using Lemma \ref{88} once more,
we obtain that
$$ \|f_s(x)\|_{L^2(\R^N)}\le
|\nu_s|\,\|u_s^p\|_{L^2(\R^N)}+
\|u_s\|_{L^2(\R^N)}\le C_3,$$
with $C_3>0$ independent of $s$ and $u_s$.
Also, from \eqref{EQU}, $\Ds u_s=f_s$, that is, recalling \eqref{FOU},
$$ |\xi|^{2s} \widehat{u_s} = \widehat{f_s}$$
and so
\begin{eqnarray*} && \|u_s\|_{2s}^2 =
\|u_s\|^2_{L^{2}(\R^N)}+\int_{\R^N} |\xi|^{4s} |\widehat{u}_s|^2\,d\xi
\le \nu_s +
\int_{\R^N} |\widehat{f}_s|^2\,d\xi\\
&&\qquad=\nu_s+\|f_s\|^2_{L^2(\R^N)}\le C_2+C_3,
\end{eqnarray*}
and the desired result plainly follows. \QED

Next result is a general approximation argument
on the fractional Laplacian:

\begin{Lemma}\label{667}
Let $s$, $\bar\s\in(0,1]$ and
\be\label{BE}
\delta>2|\bar\s-s|.\ee Then, for any $\vp\in
H^{2(\bar\s+\delta)}(\R^N)$,
$$ \| (-\Delta)^{\bar\s}\vp-(-\Delta)^{s}\vp\|_{L^2(\R^N)}
\le C_{\bar\s,\delta} |\bar\s-s|\,\|\vp\|_{2(\bar\s+\delta)} ,$$
for a suitable $C_{\bar\s,\delta}>0$.
\end{Lemma}

\proof
We start with some elementary inequalities.
First of all, if $\tau\in[0,1)$ then
$(1+\tau^{2\bar\sigma+\delta})\tau^{2|\bar\sigma-s|}\le 2\cdot 1$.
On the other hand, if $\tau\ge1$
then~$(1+\tau^{2\bar\sigma+\delta})\tau^{2|\bar\sigma-s|}\le
(2\cdot \tau^{2\bar\sigma+\delta})\tau^\delta$, thanks to~\eqref{BE}.
All in all, we obtain that, for any~$\tau\ge0$,
\be\label{67.67}
(1+\tau^{2\bar\sigma+\delta})\tau^{2|\bar\sigma-s|}\le 2
(1+\tau^{2(\bar\sigma+\delta)}).
\ee
Moreover, for any $t\in\R$,
\begin{equation}\label{7a.1}
|e^t-1|\le \sum_{k=1}^{+\infty}\frac{|t|^k}{k!}\le
\sum_{k=1}^{+\infty}\frac{|t|^k}{(k-1)!}=|t| e^{|t|}.
\end{equation}
Furthermore, the map~$(0,1)\ni\tau\mapsto \tau^{2\bar\sigma}\log\tau$
is minimized at~$\tau=e^{-1/2\bar\sigma}$ and therefore
\be\label{90.01}
|\tau^{2\bar\sigma}\log\tau|\le (2\bar\sigma e)^{-1}\qquad
{\mbox{ for any }}\tau\in(0,1).
\ee
Similarly, the map~$[1,\infty)\ni\tau\mapsto \tau^{-\delta}\log\tau$
is maximized at~$\tau=e^{1/\delta}$ and so
\be\label{90.02}
|\tau^{-\delta}\log\tau|\le (\delta e)^{-1}\qquad
{\mbox{ for any }}\tau\in[1,\infty).
\ee
By combining \eqref{90.01} and~\eqref{90.02}, we obtain that,
for any $\tau>0$,
\begin{equation}\label{7a.2}
|\tau^{2\bar\sigma} \log\tau|\le C_{\bar\sigma,\delta}\;(1+\tau^{2\bar\sigma+\delta})
\end{equation}
where
\be\label{BE3}
C_{\bar\sigma,\delta}:=
(2\bar\sigma e)^{-1}+(\delta e)^{-1}
.\ee
Thus, using \eqref{67.67},
\eqref{7a.1} and \eqref{7a.2}, we obtain that,
for any $\xi\in\R^N\setminus\{0\}$,
\begin{eqnarray}\label{LL72}\nonumber
&& ||\xi|^{2s}-|\xi|^{2\bar\s}|=|\xi|^{2\bar\s} |\xi^{2(s-\bar\s)}-1|
=|\xi^{2\bar\s}|\,|e^{2(s-\bar\s)\log|\xi|}-1|\\
&&\qquad
\le |\xi|^{2\bar\s}\,\big|{2(\bar\s-s)\log|\xi|}\big| e^{{2|\bar\s-s|\,|\log|\xi||} |}
=|\xi|^{2\bar\s}\,2|\bar\s-s|\,\big|\log|\xi|\big|\;|\xi|^{2|\bar\s-s|}
\\ &&\qquad \le
2C_{\bar\sigma,\delta} |\bar\s-s|(1+|\xi|^{2\bar\s+\delta}) \;|\xi|^{2|\bar\s-s|} \le
4C_{\bar\sigma,\delta} |\bar\s-s|(1+|\xi|^{2(\bar\s+\delta)}).
\nonumber\end{eqnarray}
As a consequence
\begin{eqnarray*} && \int_{\R^N} |[\Ds \phi-\Dsib \phi] |^2=
\int_{\R^N} ||\xi|^{2s} - |\xi|^{2\bar\s}|^2 |\widehat{\phi}|^2 \\
&&\qquad\leq {\mbox{Const}}\,C_{\bar\sigma,\delta}^2
(\bar\s-s)^2\int_{\R^N}(1+|\xi|^{4(\bar\s+\delta)})
|\widehat{\phi}|^2\,d\xi
\leq {\mbox{Const}}\, C_{\bar\sigma,\delta}\,\|\phi\|_{2(\bar\s+\delta)}^2,\end{eqnarray*}
as desired.\QED

\begin{Corollary}\label{9866:0} Fix $\sigma\in(0,1]$. Then
$\displaystyle\lim_{s\to\sigma}\nu_{s}=\nu_{\sigma}.$
\end{Corollary}

\proof
Let $s_0\in(0,1)$. Let~$s$, $s'\in(s_0,1]$, that will be taken
one close to the other, namely such that
\be\label{BE2}
s>2|s-s'|.
\ee
Let $u_s\in\calM_s$.
Since $\|u_s\|_{L^{p+1}(\R^N)}=1$,
we obtain that $ \nu_{s'}\le \|u_s\|_{s'}^2$. Hence, recalling
\eqref{BE3} and
\eqref{LL72} (used here with $\bar\sigma:=s'$ and $\delta:=s$,
and notice that \eqref{BE} is warranted by~\eqref{BE2}), we
conclude that
\begin{eqnarray*}
&& \nu_{s'}-\nu_s\le \|u_s\|_{s'}^2 -\|u_s\|_{s}^2\\
&& \qquad=\int_{\R^N}\Big(|\xi|^{2s'}-|\xi|^{2s}\Big)\,|\widehat{u_s}|^2
\le {\mbox{Const}}\,|s'-s|
\int_{\R^N}\Big(1+|\xi|^{4s}\Big)\,|\widehat{u_s}|^2\\ &&\qquad=
{\mbox{Const}}\,|s'-s|\,\|u_s\|_{2s}
\end{eqnarray*}
The constants here above only depend on the fixed~$s_0$, but not
on~$s$ and~$s'$.
Since the roles of $s$ and $s'$ may be interchanged,
and recalling Corollary~\ref{C:REG}, we obtain that
$$ |\nu_{s'}-\nu_s|\le {\mbox{Const}}\,|s'-s|$$
and the desired result plainly follows.
\QED

{F}rom now on, we will use the uniqueness and nondegeneracy results
for the local case. Namely, we recall that
there exists a unique   radial
minimizer $U_1(x)=\bar{U}_1(|x|) $ for $\nu_1$, such that
\be\label{eq:ND}
Ker(J_1''(U_1,\nu_1)) =\textrm{span}\{\de_j U_1,\,j=1,\dots,N\},
\ee
see, e.g. \cites{serrin-mcleod, mcleod, kwong}.

\begin{Lemma}\label{lem:reg} Fix $\bar{\s}\in (0,1]$. Let $s_n\in(0,1)$ be
such that $s_n\to\bar{\s}$. Let $u_{s_n}\in \calM_{s_n}$. Then there exist
$ \bar{u}\in \calM_{\bar{\s}}$ and a subsequence (still denoted by $s_n$)
such that if \be\label{E6}
\o_{s_n}(x):= u_{s_n}(x)-\bar{u},\ee we have that
$$ {\mbox{$\|\o_{s_n}\|_{2s_n}\to 0$ as $n\to \infty.$}}$$
Moreover, if
$\bar{\s}=1$ then $${\mbox{$\|\o_{s_n}\|_{2}\to 0$ as $n\to \infty$. }}$$
\end{Lemma}

\proof To alleviate the notation, we write $s$ instead of $s_n$. {F}rom
Corollary~\ref{C:REG}
we have that $ u_s$ is
bounded
in $H^t(\R^N)$ for every $t<\bar\s$. Therefore, by compactness (see
Proposition 1.1 in~\cite{lions}), we obtain that there
exists $\bar u$ such that
$$
{\mbox{$u_s\to \bar{u}$ in $L^q(\R^N)$
for every $q\in (2,2^*_{\bar\s})$.}}$$
Since we have uniform decay bounds at infinity and uniform $L^\infty$ bounds 
(recall Lemma~\ref{lem:unif}),
this and the interpolation inequality implies
that the convergence also holds for~$q\in(1,2]$, hence
\be\label{BL}
{\mbox{$u_s\to \bar{u}$ in $L^q(\R^N)$
for every $q\in (1,2^*_{\bar\s})$.}}\ee
In particular,
$\|\bar u_{\bar\s}\|_{L^{p+1}(\R^N)}=1$ and $\bar{u}$ is radially
symmetric.
What is more, by Fatou lemma, it follows that $\bar u\in H^{\bar\s}(\R^N)$ because $\int_{\R^N}|\xi|^{2s}|\widehat{u_s}|^2d\xi\leq \nu_s\leq Conts$.
Also, by \eqref{EQU}, \be\label{vb67}\int_{\R^N}u_s\Ds
\phi+\int_{\R^N}u_s\phi=\nu_s\int_{\R^N}u_s^p\phi\quad\forall \phi\in
C^\infty_c(\R^N). \ee
Using Lemma \ref{667},
$$ \int_{\R^N} |[\Ds \phi-\Dsib \phi] |^2\,d\xi\leq {\mbox{Const}}\,
(\bar\s-s)^2\int_{\R^N}(1+|\xi|^{4})
|\widehat{\phi}|^2\,d\xi
\leq {\mbox{Const}}\|\phi\|_{2}^2. $$
Hence we can pass to the limit in~\eqref{vb67} and
conclude that $\bar{u}$ is a distributional solution to the equation
\be\label{E7}\Dsib \bar u+
\bar u=\nu_{\bar\s} \bar u^p\ee that belongs to $H^{\bar\s}(\R^N)$.

So, by
testing the equation against~$u$ itself, we see
that~$\|u\|_{\bar \s}^2=\nu_{\bar\s}\|u\|^{p+1}_{L^{p+1}(\R^N)}=\nu_{\bar\s}$,
hence $\bar{u}$ is a minimizer for $ \nu_{\bar\s}$.

Furthermore, by \eqref{EQU}, \eqref{E6} and \eqref{E7},
\be\label{eq:dso}
\begin{array}{ccc}
\Ds\o_s+\o_s&=& \nu_s [(\bar{u}+\o_s)^{p} -\bar{u}^{p}] \hspace{3cm} \\
& &+[\Dsib \bar{u}- \Ds \bar{u}]+(\nu_s-\nu_{\bar\s}) \bar{u}^p.
\end{array}
\ee
Also, from the fundamental theorem of calculus
\begin{eqnarray*}
(\bar{u}+\o_s)^{p} -\bar{u}^{p} &=&
\int_{0}^{1} \frac{d}{dt}(\bar{u}+t\o_s)^{p}dt
\\ &=&
p\o_s \int_{0}^{1}(\bar{u}+t\o_s)^{p-1}dt,
\end{eqnarray*}
so that using \eqref{E6} and \eqref{eq:unif}
\be\label{E8}\begin{split}
&|(\bar{u}+\o_s)^{p} -\bar{u}^{p}  |\leq p|\o_s|(\|{u}_s\|_{L^\infty(\R^N)} +2\|\bar{u}\|_{L^\infty(\R^N)})^{p-1}\\ &\qquad\leq {\mbox{Const}}\; |\o_s|.
\end{split}\ee
Next we observe that, since $\bar{u},u_{s}\in C^2(\R^ N)$, \eqref{eq:dso} holds pointwise and thus, by \eqref{E8}, we obtain
\be\label{E9}\begin{split}
& \|\Ds\o_s\|_{L^2(\R^N)}^2\\ &\quad\leq \|\o_s\|_{L^2(\R^N)}^2+{\mbox{Const}}\;(|\bar\s-s|^2+ |\nu_{\bar{\s}}-\nu_s|^2+\|\o_s\|_{L^2(\R^N)}^2)\to 0
\end{split}\ee
as $ s\to\bar{\s}$. This and \eqref{BL} imply that
$\|\o_s\|_{2s}\to 0 $ as $s\nearrow \bar{\s}$, as desired.

Next we consider the case $\bar\s=1$.
By \eqref{LIN} and
\eqref{eq:unif}
we have that
for every $s$ close to 1
$$
\|\de_j u_{s}\|_{2s}\leq {\mbox{Const}}.
$$
{F}rom this, \eqref{E7} and \eqref{E6}, we deduce that
$$
\|\de_j \o_{s}\|_{2s}\leq {\mbox{Const}}.
$$
In particular $ \|\o_s\|_{2s+1}$ is uniformly bounded.
We let $f_s$ be the right hand side of \eqref{eq:dso} so that
$$
\Ds\o_s+\o_s=f_s
$$
and so
$$
-\D\o_s+\o_s=f_s+[ -\D\o_s- \Ds\o_s ].
$$
Using Lemma \ref{667}, we conclude that, for every $\d\in(0,1/4)$,
$$
\int_{\R^N}[ -\D\o_s- \Ds\o_s ]^2\leq C_{N,\d}(1-s)\|\o_s\|_{2+\d}\leq  \|\o_s\|_{2s+1}\leq (1-s) {\mbox{Const}},
$$
provided $s$ is close to 1.
Also, by recalling \eqref{E9} and \eqref{BL}, we obtain that
$\|f_s\|_{L^2(\R^N)}\to 0$ as $s\nearrow 1$, and therefore
$\|\o_s\|_{2}\to 0 $. \QED


\subsection{Local realization of $\Ds$ for $s\in(0,1)$}

Following~\cite{CaSi},
we recall here an extension property that provides a local realization of the
fractional Laplacian by means of a divergence operator in a higher dimension halfspace.
Namely, given $u\in H^s(\R^N)$, there exists a unique $\calH(u)\in H^1(\RNp;t^{1-2s})$ such that
\be\label{SSIL}
\begin{cases}
\div(t^{1-2s}\n \calH(u))=0 \quad\textrm{ in }\RNp,\\
\calH(u)=u  \quad\textrm{ in }\R^N,\\
\lim_{t\searrow0} t^{1-2s}\calH(u)_t:= t^{1-2s}\calH(u)_t=\k_s \Ds u
\quad \textrm{ on }\R^N,
\end{cases}
\ee
where $\k_s$ is a positive normalization constant.
Equivalently for every $\Psi\in H^1(\RNp;t^{1-2s}) $
\be
\int_{\RNp} \n \calH(u)\cdot\n\Psi\,t^{1-2s}dt\,dx=\k_s \int_{\R^N}|\xi|^{2s}\widehat{u}\widehat{\Psi}d\xi,
\ee
where here and hereafter we denote the trace of a function with the same letter.
From now on, $\calH$ will denote the  $s$-harmonic operator.
Moreover, the trace property holds, i.e. for any $\Phi\in  H^1(\RNp;t^{1-2s})$,
the trace $\Phi$ on $\R^N$
belongs to $H^s(\R^N)$. As $ \calH(tr(\Phi)):= \calH(\Phi)$ has minimal Dirichlet  energy, it follows that
$$
\int_{\R^{N+1}_+}|\n \Phi|^2t^{1-2s}dt\,dx    \geq  \int_{\R^{N+1}_+}|\n \calH(\Phi)|^2t^{1-2s}dt\,dx  = \k_s\int_{\R^N}|\xi|^{2s}|\widehat{\Phi}|^2d\xi.
$$
Hence $\calH(u_s)$ is radially symmetric with respect to the $x$  variable and it is  a minimizer for
\be\label{eq:nusH}
\nu_s=\inf_{U\in H^1(\RNp;t^{1-2s}) }
\frac{\displaystyle\k_s^{-1}  \int_{\displaystyle\R^{N+1}_+}|\n U|^2t^{1-2s}dt\,dx +  \int_{\R^{N}}| U|^2dx}{\displaystyle\left(\int_{\R^N}|U|^{p+1}dx\right)^{2/(p+1)}}
\ee
and, by \eqref{SSIL},
\be\label{eq:GroundH1s}
\begin{cases}
\div(t^{1-2s}\n \calH(u_s))=0 \quad\textrm{ in }\RNp\\
\k_s^{-1} t^{1-2s}\calH(u)_t+\calH(u) =\nu_s \calH(u)^p \quad \textrm{ on }\R^N.
\end{cases}
\ee
In this setting, we define
$$
\calJ_s(U,\nu):=\frac{1}{2} \int_{\R^{N+1}_+}|\n U |^2t^{1-2s}dt\,dx+\frac{\k_s }{2}\int_{\R^N}U^2 dx-\frac{\nu\k_s }{p+1}\int_{\R^N}|U|^p dx.
$$
%

\section{Nondegenracy}\label{S:2}

\subsection{Preliminary observations}

In this section, we assume that  $u_s
\in\calM_s$ and we prove that it is nondegenerate for $s$ sufficiently close to $1$.
For this, we denote by $\perp_s$ the orthogonality relation in $H^s(\R^N)$
and we start by estimating the second variation of the functional.

\begin{Lemma}\label{lem:Jsecnonneg}
For every $\phi\perp_s u_s$ we have that
\be\label{0966}
0\le J_s''(u_s,\nu_s)[\phi,\phi]=\|\phi\|^2_s-p\nu_s\int_{\R^N}u_s^{p-1}\phi^2\,dx.
\ee
\end{Lemma}
\proof
Let $\e>0$.
Since~$\phi\perp_s u_s$, we have
\be\label{B6767}
\|\e\phi+u_s\|^2_s=\e^2\|\phi\|^2_s+\|u_s\|^2_s .\ee
Also, by a Taylor expansion we obtain
\be\label{TA}\begin{split}
&\int_{\R^N}|\e\phi+u_s|^{p+1}
\\&\qquad= \int_{\R^N}|u_s|^{p+1}+\e (p+1)
\int_{\R^N}u_s^{p}\phi+\frac{\e^2p(p+1)}{2}\int_{\R^N}u_s^{p-1}\phi^2+O(\e^3).
\end{split}\ee
Furthermore, by testing~\eqref{EQU} against~$\phi$ and
using again that~$\phi\perp_s u_s$, we conclude that
$$\int_{\R^N}u_s^{p}\phi=0,$$ hence the first order in~$\e$ in~\eqref{TA} vanishes.
Consequently, recalling also that functions in~$\calM_s$ are normalized with~$\|u\|_{L^{p+1}(\R^N)}=1$,
we write~\eqref{TA} as
\be\label{TA2}\int_{\R^N}|\e\phi+u_s|^{p+1}
= 1+\frac{\e^2p(p+1)}{2}\int_{\R^N}u_s^{p-1}\phi^2+O(\e^3). \ee
Now we recall the Taylor expansion
\be\label{TA3}
\frac{1}{(1+x)^{2/(p+1)}}=1-\frac{2}{p+1}x+O(x^2)
\ee
for small~$x$. Thus, by inserting~\eqref{TA2} into~\eqref{TA3}, we obtain
$$ \frac{1}{\left(\displaystyle \int_{\R^N}|\e\phi+u_s|^{p+1}\right)^{2/(p+1)}}
=1-{\e^2p}\int_{\R^N}u_s^{p-1}\phi^2+O(\e^3).$$
{F}rom this and~\eqref{B6767} we obtain
\begin{eqnarray*}
&& \frac{ \displaystyle \|\e\phi+u_s\|^2_s   }{\left(\displaystyle \int_{\R^N}|\e\phi+u_s|^{p+1}\right)^{2/(p+1)}}
\\&&\qquad=\Big(1-{\e^2p}\int_{\R^N}u_s^{p-1}\phi^2+O(\e^3)\Big)\,\Big( \e^2\|\phi\|^2_s+\|u_s\|^2_s\Big)
\\&&\qquad=\|u_s\|^2_s+\e^2 \Big( \|\phi\|^2_s-p\|u_s\|^2_s \int_{\R^N}u_s^{p-1}\phi^2\Big)+O(\e^3).\end{eqnarray*}
Then the desired result follows since~$u_s$ attains the minimal value~$\nu_s=\|u_s\|^2_s$.
\QED


\begin{Lemma}\label{lem:ext-pos}
Let $\Phi\in  H^1(\RNp;t^{1-2s}) $ be such that
 \be\label{eq:perpH}
 \k_s^{-1}\int_{\R^{N+1}_+}\n \Phi\cdot \n \calH(u_s) t^{1-2s}dt\,dx+\int_{\R^{N}}\Phi\calH(u_s)dx =0.
\ee
Then
\be\label{eq:Phi}
\calJ''_s(\calH(u_s))[\Phi,\Phi]= \k_s^{-1} \int_{\R^{N+1}_+}|\n \Phi|^2t^{1-2s}dz+\int_{\R^{N}}\Phi^2dx-p\nu_s \int_{\R^N}u_s^{p-1}\Phi^2dx\geq 0.
\ee
In particular for any $g\in  H^1(\R^2_{++};t^{1-2s}r^{N-1})$
\begin{eqnarray}\label{90:98}
\nonumber
A_1(g,g)&:=& \int_{\R^2_+}  g_t^2t^{1-2s}r^{N-1}dtdr +  \int_{\R^2_{++}}  g_r ^2t^{1-2s}r^{N-1}dtdr\\
&&+ (N-1)  \int_{\R^2_{++}}  g^2t^{1-2s}r^{N-3}dtdr
+\k_s   \int_{\R_+}   g^2r^{N-1}dr\\
&&-p\nu_s \k_s  \int_{\R_+}  u_s^{p-1} g^2r^{N-1}dr\geq0.\nonumber
\end{eqnarray}
\end{Lemma}
\proof
The proof of \eqref{eq:Phi} is similar to
the proof of {Lemma} \ref{lem:Jsecnonneg}, since $\calH(u_s)$ minimizes \eqref{eq:nusH}.
Next, let $g\in  H^1(\R^2_{++};t^{1-2s}r^{N-1})$ and define $\Phi(x):=g(t,|x|)\frac{x^i}{|x|} $.
Since $\calH(u_s)$ is radial in the $x$ variable,
$\Phi $ satisfies \eqref{eq:perpH} by odd symmetry. Then
\eqref{eq:perpH}, \eqref{eq:Phi} and the use of polar coordinates yield~\eqref{90:98}.
\QED

\begin{Lemma}\label{lem:dimKer}
Let $w\in Ker J_s''(u_s,\nu_s)$. Then   $$w=w_0(|x|)+ \sum_{i=0}^Nc^i\de_iu_s,$$
where $$w_0(r)=\int_{S^{N-1}}w(r\th)d\s(\th)$$ and $c^i\in\R$.
\end{Lemma}

\proof
Let $w\in Ker(J_s''(u_s,\nu_s))$ which means
$$
\Ds w+w-p\nu_s u_s^{p-1}w=0\quad \textrm{ in }\R^N.
$$
Let $\calH(w)\in   H^1(\RNp;t^{1-2s})$ be the $s$-harmonic extension of $w$ which  satisfies
\be\label{eq:calHvweak}
\k_s^{-1}\int_{\RNp} \n \calH(w)\cdot\n\Psi\,t^{1-2s}dt\,dx+\int_{\R^N}\calH(w)\Psi \,dx - p\nu_s \int_{\R^N}u_s^{p-1}  \calH(w)\Psi dx=0,
\ee
for all $\Psi\in  H^1(\RNp;t^{1-2s})$.
Now we consider the spherical harmonics on $\R^N$ for $N\geq 2$, i.e.
the solution of the classical eigenvalue problem
$$-\D_{S^{N-1}} Y^i_k=\l_k Y^i_k \ {\mbox{ on $S^{N-1}$.}}$$
We let $n_k$ be the multiplicity of $\l_k$. It is known that $n_0=1$ and $n_1=N$
(see e.g. formulae~(3.1.11) and~(3.1.12) in~\cite{SPH}). In addition
$\l_0=0$, $\l_1=N-1$ and $\l_k> N-1$ for $k\geq 2$. Also $Y_0$ is constant, while
$$ Y^i_1=\frac{x^i}{|x|}\ {\mbox{ for }} \ i=1,\dots, N.$$
With this setting, we decompose $\calH(w)$ in the spherical harmonics and we obtain
\be\label{1977}
\calH(w)(t,x)=\sum_{k\in \N}\sum_{i=1}^{n_k}f^k_i(t,|x|)Y_k^i\left(\frac{x}{|x|}\right),
\ee
where $f^k_i\in H^1(\R^2_+;t^{1-2s}r^{N-1})  $. By testing \eqref{eq:calHvweak}
against the function $\Psi= h(t,|x|) Y_k^i$ and using polar coordinates, we
obtain that,
for any $ h \in H^1(\R^2_+;t^{1-2s}r^{N-1}) $, any $k\in\N$ and any $i\in[1,n_k]$,
\begin{eqnarray*}
A_k(f^k_i,h)&:=& \int_{\R^2_{++}}  (f^k_i)_t h_t t^{1-2s}r^{N-1}dtdr +  \int_{\R^2_{++}}  (f^k_i)_r h_rt^{1-2s}r^{N-1}dtdr\\
&&+ \l_k  \int_{\R^2_{++}}  f^k_ih t^{1-2s}r^{N-3}dtdr
+\k_s    \int_{\R_+} f^k_i h r^{N-1}dr\\
&&-p\nu_s \k_s  \int_{\R_+}  u_s^{p-1} f^k_i h r^{N-1}dr =0.
\end{eqnarray*}
Now we observe that
$$
A_k(f^k_i,f^k_i)=A_1(f^k_i,f^k_i)+(\l_k-(N-1))\int_{\R^2_{++}}\int_{S^{N-1}} (f^k_i)^2 t^{1-2s}r^{N-3}dtdr.
$$
By   Lemma \ref{lem:ext-pos} and the fact that  $ \l_k>N-1$ for $k\geq 2$, we obtain from the identities above that
\begin{eqnarray*} && 0=A_k(f^k_i,f^k_i)=
A_1(f^k_i,f^k_i)+(\l_k-(N-1))\int_{\R^2_{++}}\int_{S^{N-1}} (f^k_i)^2 t^{1-2s}r^{N-3}dtdr\\ &&\geq
(\l_k-(N-1))\int_{\R^2_{++}}\int_{S^{N-1}} (f^k_i)^2 t^{1-2s}r^{N-3}dtdr\ge0.
\end{eqnarray*}
As a consequence, $f^k_i=0$ for every $k\geq 2$.
Accordingly, \eqref{1977} becomes
$$ \calH(w)(t,x)=\sum_{i=1}^{N}f^1_i(t,|x|)Y_k^i\left(\frac{x}{|x|}\right).$$
To complete the proof we need to characterize~$f^1_i$.
For this, we notice that, for $i=1,\dots,N$,  the function
$$f^1_i(t,r)=\int_{S^{N-1}}\calH(w)(t,r\th)\th^i d\s(\th)$$ satisfies
 $f^1_i(t,0)=0$ and
 \begin{eqnarray}\label{9878644}\nonumber
A_1(f^1_i,h)&=& \int_{\R^2_{++}}  (f^1_i)_t h_t t^{1-2s}r^{N-1}dtdr +  \int_{\R^2_{++}}  (f^1_i)_r h_rt^{1-2s}r^{N-1}dtdr\\
&&+(N-1)  \int_{\R^2_{++}}  f^1_ih t^{1-2s}r^{N-3}dtdr
+\k_s     \int_{\R_{++}} f^1_i h r^{N-1}dr\\
&&-p\nu_s \k_s  \int_{\R_{++}}  u_s^{p-1} f^1_i h r^{N-1}dr =0,\nonumber
\end{eqnarray}
for every $h\in H^1(\R^2_+;t^{1-2s}r^{N-1})$, due to~\eqref{eq:calHvweak}.

Now we define $\bar{U}(t,|x|)=\calH(u_s)(t,x) $. Then we have
$$
\begin{cases}
 \div(t^{1-2s}r^{N-1} \n \bar{U} )=0 \quad \textrm{ in }\R^2_{++}\\
\,\\
\lim_{t\searrow0}
- t^{1-2s}r^{N-1}\bar{U}_t+ \k_sr^{N-1}\bar{U}= \k_sr^{N-1}\bar{U}^p
\quad\textrm{ on } \R_{+}.\\ \, \\
\lim_{r\searrow0} r^{N-1}  \bar{U}_r(t,0)=0.
\end{cases}
$$
We set $V:=\bar{U}_r$ and we differentiating the above equation with respect to $r$. We obtain
\be\label{eq:eqforV}
\begin {cases}
 -\div(t^{1-2s}r^{N-1} \n V )+(N-1)t^{1-2s}r^{N-3}V=0 \quad\textrm{ in }
\R^2_{++}\\ \, \\
\lim_{t\searrow0}- t^{1-2s}r^{N-1}V_t+ \k_sr^{N-1}V=
\k_spr^{N-1}\bar{U}^{p-1} V\quad\textrm{ on } \R_+\\
\, \\\lim_{r\searrow0} r^{N-1}  V(t,0)=0.
\end {cases}
\ee
Since $\bar{U}_r$ does not change sign, we may assume that $V<0$ on $\R^2_{++}$.

Given $g\in C^\infty_c({\R_{++}^2}\cup\{t=0\})$, we   define $$\psi:= \frac{g}{V}\in H^1(\R^2_{++};t^{1-2s}r^{N-1}).$$
Simple computations  show that
$$
|\n g|^2=|V \n\psi|^2+\n V \cdot \n({V }\psi^2).
$$
Hence we have
\begin{eqnarray*}&&
\int_{\R^2_+} |\n g|^2t^{1-2s}r^{N-1}dtdr
\\ &&\quad=\int_{\R^2_+}|V \n\psi|^2 t^{1-2s}r^{N-1}  dtdr+\int_{\R^2_{++}} \n({V }\psi^2)\cdot(t^{1-2s}r^{N-1} \n V)dtdr.
\end{eqnarray*}
Integrating by parts, by  using the above identities and \eqref{eq:eqforV}, we get
\begin{eqnarray*}
\int_{\R^2_{++}} |\n g|^2t^{1-2s}r^{N-1}dtdr+(N-1)\int_{\R^2_{++}}   g^2 t^{1-2s}r^{N-3}dtdr + \k_s \int_{\R_+}  g^2r^{N-1}dr   \\
 -  \k_s p \int_{\R_{+}} u_s^{p-1} g^2r^{N-1}dr= \int_{\R^2_{++}}\left|V \n\psi\right|^2 t^{1-2s}r^{N-1}  dtdr.
\end{eqnarray*}
In particular, by density and recalling~\eqref{9878644}, we have that, for every $i=1,\dots,N$,
$$
A_1(f^1_i,f^1_i)=0\geq\int_{\R^2_{++}}\left|V \n\left( f^1_i V^{-1}\right)\right|^2 t^{1-2s}r^{N-1}  dtdr.
$$
This implies that the last term vanishes and
therefore  $$\frac{ f^1_i  }{V} \equiv c^i$$ for some constant $c^i\in\R$. We then conclude that  $ f^1_i(0,|x|)  =c^i \bar{U}_s'(|x|)$ for all $x\in\R^N$.\\

Thus, we have proved that for any $w\in Ker(J_s''(u_s,\nu_s))$
$$
\calH(w)(0,x)=w(x)=f^0_1(0,|x|)+\sum_{i=1}^{N}f^1_i(0,|x|) \frac{x^i}{|x|}=
f^0_1(0,|x|)+\sum_{i=0}^{1}c^i\de_ku_s(x),
$$
as desired.
\QED

Now we are ready to prove our nondegeneracy result for $s$ close to $1$.

\subsection{Completion of the proof of
Theorem \ref{th:nondegen}}

Let $ v_s\in Ker(J''(u_s,\nu_s))$ be a radial function.  \\

\noindent
\textbf{Claim:} If $s$ is close to $1$, we have $v_s\equiv0$.\\
Assume by contradiction that there exists a sequence $s_n$ -- still denoted by $s$ -- with $s \nearrow 1$ and
such that
${v}_s\neq 0$. Up to normalization, we can assume that $\|v_s\|_{L^{p+1}(\R^N)}=1$. By Corollary \ref{9866:0}, we know
that
$\nu_s\to \nu_1$ and $u_s\to U_1(\cdot-a)$ in $L^{p+1}$, for some $a\in \R^N$.
Since $u_s$ is symmetric with respect to the origin, $a=0$.
By H\"older inequality
\be\label{L72}
\int_{\R^N}|\xi|^{2s}|\widehat{v_s}|^2d\xi\leq \|v_s\|^2_s\leq p\nu_s\|u_s\|_{L^{p+1} }^{p-1}\|v_s\|_{L^{p+1}(\R^N)}^2=  p\nu_s\leq {\mbox{Const.}},
\ee
by Lemma \ref{88}.
Since $v_s$ is a radial sequence and bounded in $H^t(\R^N)$ for every $t\in(0,1)$, by compactness (see \cite{lions})
$ v_s\to v$ in $L^q(\R^N)$ for every $q\in (2,2^*_1)$,
and then also for~$q=2$ (by repeating the argument above~\eqref{BL}). In particular~$\|v\|_{L^{p+1}(\R^N)}=1$.
Next we observe that $v_s$
is a solution of the linearized equation
and therefore for any $\phi\in C^\infty_c(\R^N)$
$$
\int_{\R^N }v_s \Ds\phi  +\int_{\R^N }v_s \phi -p\nu_s\int_{\R^N }  u_s^{p-1}v_s \phi =0
$$
so by \eqref{eq:unif} and the fact that $\Ds\phi\to -\D \phi$ in $L^2(\R^N)$
thanks to Lemma \ref{667}, we infer that
$$
\int_{\R^N }v (-\D)\phi  +\int_{\R^N }v \phi -p\nu_1\int_{\R^N }  U_1^{p-1}v \phi =0.
$$
Applying Fatou lemma to \eqref{L72}, we get $v\in H^1(\R^N)$. We then conclude that $v$ is radial, nontrivial
and belongs to  $Ker(J''(U_1,\nu_1))$. This is  clearly a contradiction and the claim is proved. \QED

%

\section{Uniqueness (preliminary observations)}\label{S:3}

\subsection{Preliminary observations}
Now we prove Theorem~\ref{th:Uniq}.
The first part of the proof of the following result is quite standard but
the last part requires a more delicate analysis on radial functions.

\begin{Lemma}\label{lem:inv}
Let~$\L_s :=(Ker(J_s''(u_s,\nu_s))\oplus\R u_s)^{\perp_s} $.
\begin{enumerate}
\item We have
\be\label{eq:Jssec} J_s''(u_s,\nu_s)[u_s,u_s] =(1-p)\|u_s\|_s^2.\ee
\item There exists $s_0\in(0,1)$ such that
 for every $s\in({s_0},1)$ and every minimizer $u_s$ for $\nu_s$
\be\label{II}
K(s,u_s):=\inf_{\phi\in \L_s\setminus\{0\}}\frac{ J_s''(u_s,\nu_s)[\phi,\phi]   }{\|\phi\|_s^2}>0.
\ee
\item
Let
$$
\L_s^r:=\{ \phi\in H^s_{rad}(\R^N), \phi\perp_s u_s   \}
$$
and
$$
K_r(s,u_s):=\inf_{\L_s^r\setminus\{0\}}\frac{ J_s''(u_s,\nu_s)[\phi,\phi]   }{\|\phi\|_s^2}.
$$
Then there exits $s_0\in(0,1)$ such that
\be\label{ZA}
\inf_{s\in(s_0,1]}\,\inf_{u\in \calM_{s} } K_r(s,u)>0.
\ee
\end{enumerate}
\end{Lemma}
\proof
The statement in~\eqref{eq:Jssec}
is immediate from~\eqref{0966}.

Now we prove~\eqref{II}.
We first show that for any $\phi \in \L_s$
\begin{equation}\label{98f8}
J_s''(u_s,\nu_s)[\phi,\phi]=0 \Longrightarrow \phi\equiv 0.
\end{equation}
That is to say  that $J_s''(u_s,\nu_s) $ defines a scalar product on $\L_s$ by Lemma \ref{lem:Jsecnonneg}.
For this,
assume that $\phi\in \L_s$ and
$$
J_s''(u_s,\nu_s)[\phi,\phi]=0.
$$
Pick $\psi \in H^s(\R^N)$ such that $\psi\perp_s u_s$. By Lemma \ref{lem:Jsecnonneg}
$$
J_s''(u_s,\nu_s)[\phi+\e \psi ,\phi+\e\psi ]\geq 0.
$$
Hence
\begin{eqnarray*}
0\le J_s''(u_s,\nu_s)[\phi,\phi]+ 2\e J_s''(u_s,\nu_s)[\phi,\psi]+ \e^2  J_s''(u_s,\nu_s)[ \psi ,\psi ]\\
\hspace{3cm}=  2\e J_s''(u_s,\nu_s)[\phi,\psi]+ \e^2  J_s''(u_s,\nu_s)[ \psi ,\psi ].
\end{eqnarray*}
Then we conclude that \be\label{AB}{\mbox{$J_s''(u_s,\nu_s)[\phi,\psi]=0$ for any $\psi\perp_s u_s$.}}\ee
Now we observe that, since $\phi\perp_s u_s $, we deduce from \eqref{EQU} that
$$ 0=\la \phi,u_s\ra_s=\nu_s\int_{\R^N}u_s^{p}\phi$$
and so
\begin{eqnarray*}
J_s''(u_s,\nu_s)[\phi,u_s]=\la \phi,u_s\ra_s-p\nu_s\int_{\R^N}u_s^{p}\phi=0.
\end{eqnarray*}
This and \eqref{AB} yield that
$\phi \in Ker(J_s''(u_s,\nu_s)) $. Since also~$ \phi\perp_s  Ker(J''_s(u_s,\nu_s))$ it follows that $\phi=0 $, and~\eqref{98f8}
is proved.

Now we end the proof of statement 2.
Assume by contradiction that there exits  a sequence $\phi_n\in \L_s$  such that $\|\phi_n\|_s=1$ and
\be\label{j.88}
J_s''(u_s,\nu_s)[\phi_n,\phi_n]\to 0 \textrm{  as } n\to \infty.
\ee
Let~$\phi$ be the weak limit of $\phi_n$ in $H^s(\R^N)$. Then, by Lemma \ref{lem:Jsecnonneg}, we have that
$$
0\leq J_s''(u_s,\nu_s)[\phi,\phi]\leq \liminf  J_s''(u_s,\nu_s)[\phi_n,\phi_n]=0.
$$
We deduce from this and~\eqref{98f8}
that $\phi=0$, that is
\be\label{78:j:L} {\mbox{$\phi_n$ converges to $0$ weakly in $H^s(\R^N)$.}}\ee
Now,
since $u_s^{p-1}\in L^{\frac{p+1}{p-1}} (\R^N)$, given $\e>0$  there exists $w_\e\in C^\infty_c(\R^N)$ such that
\be\label{eq:apprius}
\|u_s^{p-1}-w_\e\|_{  L^{\frac{p+1}{p-1}} (\R^N)}<\e.
\ee
Now we use \eqref{78:j:L} and the compactness results
in fractional Sobolev spaces (see, e.g., Theorem 7.1
in \cite{guida}): we obtain that $\phi_n$ converges to $0$
in $L^2_{loc}(\R^n)$ and therefore
\be\label{78:89:99} \left|
\int_{\R^N}w_\e\phi_n^2\right|\le\|w_\e\|_{L^\infty(\R^N)}\,\|\phi_n\|^2_{L^2(\textrm{Supp}\,{w_\e})}
\,\to\,0\ee
as $n\to\infty$.
Also, by H\"older inequality
\begin{eqnarray*}
\left|\int_{\R^N}u_s^{p-1}\phi_n^2\right|&\leq& \|u_s^{p-1}-w_\e\|_{  L^{\frac{p+1}{p-1}} (\R^N)}^{p-1}\|\phi_n\|_{L^{p+1}(\R^N)}^2+\int_{\R^N}w_\e\phi_n^2\\
&\leq &\e^{p-1}\nu_s^{-1}+
\int_{\textrm{Supp}\,{w_\e}}w_\e\phi_n^2.
\end{eqnarray*}
This, \eqref{eq:apprius} and
\eqref{78:89:99}
imply that $$ \int_{\R^N}u_s^{p-1}\phi_n^2=o(1)
\ {\mbox{ as }}\ n\to\infty.
$$ Hence, recalling Lemma \ref{88}, we obtain
\begin{eqnarray*}
J_s''(u_s,\nu_s)[\phi_n,\phi_n]= \|\phi_n\|_s^2-p \nu_s\int_{\R^N}u_s^{p-1}\phi_n^2=1+ o(1).
\end{eqnarray*}
But this is in contradiction with \eqref{j.88} and the proof of~\eqref{II}
is complete.

Now we prove~\eqref{ZA}.
Assume by contradiction that for every $s_0\in(0,1)$
$$
\inf_{s\in(s_0,1]}\,\inf_{u_s\in\calM_s}
K_r(s,u_s)=0.
$$
Then there exits a sequence $s_n \nearrow 1 $ and radial
minimizers $u_{s_n}$ for $\nu_{s_n}$ such that
\be\label{eq:Ksusnto0}
K_r(s_n,u_{s_n})\to 0 \textrm{  as } n\to\infty.
\ee
Up to a
subsequence, and recalling
Corollary \ref{9866:0},
we may assume that $\nu_{s_n}\to \nu_1$ and,
by Lemma \ref{lem:reg}, that
\be\label{eq:usnoU1H2}
\|u_{s_n}-U_1\|_{2} \to 0 \quad \textrm{ as } n\to\infty.
\ee
For fixed $n\in\N$, by the Eckeland variational principle (see~\cite{eke}) together
with the Riesz representation theorem, we obtain that
there exist $f_{n,m}\in \L_s^r$ and  a minimizing sequence  $\psi_{n,m}\in \L_s^r$
for $K_r(s_n,u_{s_n})$ such that
$$
\|\psi_{n,m}\|_{s_n}=1,\quad\forall m\in\N
$$
and
\be\label{eq:Eckland}
J_s''(u_{s_n},\nu_{s_n})[\psi_{n,m},v]-K_r(s_n,u_{s_n}) \la\psi_{n,m} ,v \ra_{s_n}=\la f_{n,m},v \ra_{s_n}, \quad\forall v\in \L_{s_n}^r,
\ee
where $\|f_{n,m}\|_{s_n}\to 0$ as $m\to \infty$. Then there exists a sequence of
sub-indices  $m_n$ such that $\|f_{n,m_n}\|_{s_n}\to 0$ as $n\to\infty$. In particular, from~\eqref{eq:Eckland} we have
\be\label{eq:Ecklandpsinmn}
J_s''(u_{s_n},\nu_{s_n})[\psi_{n,m_n},v]-K_r(s_n,u_{s_n}) \la\psi_{n,m_n} ,v \ra_{s_n}=\la f_{n,m_n},v \ra_{s_n}.
\ee
Let $w\in C^\infty_c(\R^N)\cap \L_1^r$. Then,
from \eqref{LL72}
and \eqref{eq:usnoU1H2} we have $\la w,u_{s_n} \ra_{s_n}=o(1)\|w\|_{2+1}$
and
$$
\int_{\R^N}(1+|\xi|^{2s_n})(\imath\xi^j)u_{s_n}w=o(1)\|w\|_{2+2}\quad\forall j=1,\dots,N.
$$
We define
$$
v_n=w-\frac{ \la w,u_{s_n} \ra_{s_n} }{\|u_{s_n}\|_{s_n}^2 }u_{s_n}.
$$
By construction $v_n\in \L_{s_n}^r$. Using it as test function in \eqref{eq:Ecklandpsinmn} and recalling that
$\psi_{n,m_n}\in \L_{s_n}^r$, we get
\be\label{eq:winprod}
J_{s_n}''(u_{s_n},\nu_{s_n})[\psi_{n,m_n},w]-K_r(s_n,u_{s_n}) \la\psi_{n,m_n} ,w \ra_{s_n}=o(1).
\ee
Since $\|\psi_{n,m_n}\|_{s_n}=1$, we may assume that up to a subsequence $\psi_{n,m_n} \rightharpoonup \psi$
in $H^t(\R^N)$ for every fixed $t\in(0,1)$.
Passing to the limit in \eqref{eq:winprod}
and recalling \eqref{eq:Ksusnto0}, we get
$$
J_1''(U_1,\nu_1)[\psi,w] =0 \qquad\quad\forall w\in C^\infty_c(\R^N)\cap \L_1^r.
$$
Since, by Fatou's lemma,
$\psi\in H^1(\R^N)$, the latter identity implies that $\psi=0$, because the case $s=1$ is nondegenerate
and~$\psi\in \L_1^r$.

That is, $\psi_{n,m_n} \rightharpoonup \psi=0$
in $H^t(\R^N)$ for every fixed $t\in(0,1)$ and so, by compactness,
$${\mbox{$ \psi_{n,m_n}\to 0$ in $L^{p+1}(\R^N)$.}}$$
Also,
by \eqref{eq:Eckland}, we have
$$
J_s''(u_{s_n},\nu_{s_n})[\psi_{n,m_n},\psi_{n,m_n}]-K_r(s_n,u_{s_n}) \|\psi_{n,m_n} \|_{s_n}^2=o(1)
$$
and, by H\"older inequality
\begin{eqnarray}\label{777}\nonumber
&& \left|\int_{\R^N}u_{s_n}^{p-1}\psi_{n,m_n}^2 \right|\le\left(
\int_{\R^N}u_{s_n}^{p+1}\right)^{\frac{p-1}{p+1}}\left(
\int_{\R^N}\psi_{n,m_n}^{p+1}\right)^{\frac2{p+1}}\\&&\qquad= \|\psi_{n,m_n}\|_{L^{p+1}(\R^N)}=o(1)
\end{eqnarray}
as $n\rightarrow+\infty$. Therefore
$$ 1-p\int_{\R^N}u_{s_n}^{p-1}\psi_{n,m_n}^2
=\|\psi_{n,m_n} \|_{s_n}^2-p\int_{\R^N}u_{s_n}^{p-1}\psi_{n,m_n}^2=o(1).
$$
Hence, passing to the limit and using~\eqref{777},
we get $1-0=0$, that is a contradiction.~\QED

\section{Uniqueness (construction of pseudo-minimizers and completion
of the proof)}\label{S:4}

\subsection{Construction of pseudo-minimizers}
Pick $u_s$ a radially symmetric  even minimizer for $\nu_s$.
Define the mapping
\be
\Phi_{s}: H^s_{rad}(\R^N)\to H^{s}_{rad}(\R^N)
\ee
by
\be\label{eq:defPhias}
\Phi_{s}(\o)= J_s'\left( U_1+\o,\nu_s \right).
\ee

As customary,
by \eqref{eq:defPhias}, we mean: for all $w\in H^s_{rad}(\R^N)$
\be
\la \Phi_{s}(\o),w\ra_{s}=  J_s'\left( U_1+\o,\nu_s \right)[ w].
\ee
\begin{Lemma}\label{lem:invPHI}
 For every $f\in H^s_{rad}(\R^N)$, there exists a unique $\bar{w}^{s}\in H^s_{rad}(\R^N)$ such that
\be\label{51A}
 \la \Phi_{s}'(0)[\bar{w}^{s}],w\ra_{s}=\la f,w \ra_s\quad\forall w\in H^s_{rad}(\R^N).
\ee
In addition there exists a constant $C_1>0$ such that
\be\label{51B}
\|(\Phi_{s}'(0) )^{-1}\|\leq C_1\quad\forall s\in(s_0,1)        .
\ee
\end{Lemma}
\proof
We observe that
\begin{eqnarray*}
 \la \Phi_{s}'(0)[w'],w\ra_s=  (J_s''\left( U_1,\nu_s \right)[w', w]).
\end{eqnarray*}
Hence solving the equation
$$
 \la \Phi_{s}'(0)[\bar{w}],w\ra_{s}=\la f,w \ra_s\quad\forall w\in H^s_{rad}(\R^N)
$$
is equivalent  to find a solution~$\bar{w}$ to the equation
\begin{eqnarray}\label{eq:invsolv}
 J_s''\left( U_1,\nu_s \right)[\bar{w}, w]=\la f,w \ra_s,
\end{eqnarray}
for any~$w\in H^s_{rad}(\R^N)$.
To this scope, we observe that, for every $w\in H^s_{rad}(\R^N)$,
\begin{eqnarray}\label{P23}\nonumber
|(J_s''(U_1,\nu_s)-J_s''(u_{s},\nu_s))[w,w]|&=&\nu_sp\left|\int_{\R^N}(u_{s}^{p-1}-U_1^{p-1})w^2\right|\\
&\leq&\nu_sp\| u_s^{p-1}-U_1^ {p-1}\|_{L^{\frac{p+1}{p-1}} (\R^N)}\|w\|_{L^{p+1}(\R^N)}^2.
\end{eqnarray}
{F}rom Lemma \ref{lem:reg} and Corollary \ref{9866:0}
we know that $\|u_s-U_1\|_s\to 0$ and $\nu_s\to\nu_1$ as $s\nearrow 1$ . This implies that $u_s\to U_1$ in $ L^{p+1}(\R^N)$ and thus 
 we have
$$
 u_s^{p-1}\to U_1^ {p-1}\quad \textrm{ in } L^{\frac{p+1}{p-1}} (\R^N).
$$
Therefore, from~\eqref{P23},
\be\label{eq:JUsclous}
|(J_s''(U_1,\nu_s)-J_s''(u_{s},\nu_s))[w,w]|=o(1)\|w\|_{L^{p+1}(\R^N)}^2 .
\ee
This together with~\eqref{ZA} and \eqref{eq:Jssec} in Lemma \ref{lem:inv}  implies that there exist $C,s_0>0$ such that for all $s\in (s_0,1)$
\be\label{eq:PEP1}
| J_s''(U_1,\nu_s)[v,v]|\geq C \|v\|_s^2\quad \forall v\in H^s_{rad}(\R^N) .
\ee
Hence, by the Lax-Milgram theorem, there exits a unique $ \bar{w}^{s}\in H^s_{rad}(\R^N)$ such that
$$ J_s''(U_1,\nu_s)[\bar{w}^{s}]= f
$$
and by \eqref{eq:PEP1}
$$\|\bar{w}^{s}\|_s\leq C \|f\|_s,$$
which gives the desired result.
\QED
%
\begin{Proposition}\label{prop:finite}
For $r>0$ and $s>0$, we set
$$\calB_{r,s}=\Big\{\o\in H^s_{rad}
(\R^N)\,:\,\|\o\|_s\leq r \max\{ 1-s,|\nu_1-\nu_s|\}\Big\}.
$$
Then there exist $s_0\in(0,1)$, $r_0>0$ such that for any $s\in(s_0,1 )$, there exists a unique function
$\o^{s}\in \calB_{r_0,{s_0}}$
such that
$$
\Phi_{s}(\o^{s})=0 .
$$
\end{Proposition}
\proof
We transform the  equation $ \Phi_{s}(\o)=0$ to a  fixed point equation:
\be\label{eq:fixed}
\o=-(\Phi_{s}'(0))^{-1}\left\{\Phi_{s}(0)+Q_{s}(\o) \right\},
\ee
where
$$
Q_{s}(\o):= \Phi_{s}(\o)-\Phi_{s}(0) -\Phi_{s}'(0)[\o].
$$
Notice that the definition above is
well-posed thanks to~\eqref{51B}.
We observe
that if $\o\in H^s_{rad}(\R^N)$ then the mapping  $\o\mapsto (\Phi_{s}'(0))^{-1}\left\{\Phi_{s}(0)+Q_{s}(\o) \right\} $ is
radial too, since~$U_1$ is radial.

For very $\bar{\o}\in H^s_{rad}(\R^N)$, we set
\begin{eqnarray*}
\calN_{s}(\o)[\bar{\o}]&:=&  J_s'(U_1+\o,\nu_s )[\bar{\o}]-J_s'(U_1  ,\nu_s )[\bar{\o}]-J_s''(U_1,\nu_s )[\o,\bar{\o}]\\
&=&\nu_s\left(-\int_{\R^N}|U_1+\o|^p\bar{\o} dx+\int_{\R^N}U_1^p\bar{\o} dx  +\int_{\R^N}U_1^{p-1}\o\bar{\o} dx \right).
\end{eqnarray*}
Notice that
\be\label{7BR8}
Q_{s}(\o)= \calN_{s}(\o).\ee
Also, referring to page~128 in~\cite{AM}, we obtain
$$
|\calN_{s}(\o)[\bar{\o}]|\leq C(\|\o\|^2_s+\|\o\|^p_s)\|\bar{\o}\|_s
$$
and
$$
\|\calN_{s}(\o_1)-\calN_{s}(\o_2)\|\leq C(\|\o_1\|_s+\|\o_1\|^{p-1}_s+\|\o_2\|_s+\|\o_2\|^{p-1}_s)\|\o_1-\o_2\|_s .
$$
This, together with  ~\eqref{7BR8},
implies that  for every $\|\o_1\|_s, \|\o_2\|_s <1$,
\be\label{eq:quadr}
\|Q_{s}(\o_1)\|\leq C_3\|\o_1\|_s^{\min(2,p)}
\ee
and
\be\label{eq:contrac}
\|Q_{s}(\o_1)-Q_{s}(\o_1)\|\leq C_3\|\o_1-\o_2\|_s,
\ee
where $C_3$ is independent on $s\in (s_0,1)$.

Now we claim that
there exists a  constant  $C_2>0$
independent on  $s\in(s_0,1)$ such that
\be\label{eq:nrjsma}
\|\Phi_{s}(0)\|\leq C_2 \max\{ 1-s,|\nu_1-\nu_s|).
\ee
By~\eqref{LL72} we conclude that
$$
|J_s'(U_1 ,\nu_s)[v]- J_1'(U_1 ,\nu_1)[v]|\leq (1-s)C_{\d,N} \| U_1\|_{2-s+\d} \|v\|_s+|\nu_1-\nu_s| \|v\|_s.
$$
Since, from~\eqref{EQU}, $J_1'(U_1,\nu_1)=0 $,
we get \eqref{eq:nrjsma}.

Now we finish the proof of Proposition~\ref{prop:finite}.
We shall solve the fixed point equation \eqref{eq:fixed} in a ball of the form
$$
\calB_{r,s}=\{\o\in H^s_{rad}(\R^N)\,:\,\|\o\|_s\leq r\a_s\},
$$
where $\a_s= \max\{ 1-s,|\nu_1-\nu_s|\}$ and   $r>0$ will be fixed  in a minute.
Indeed for  $\o\in \calB_{r,s}$, we exploit~\eqref{51B}, \eqref{eq:nrjsma} and~\eqref{eq:quadr} to deduce that
$$
\|(\Phi_{s}'(0))^{-1}\left\{\Phi_{s}(0)+Q_{s}(\o) \right\}\|_s\leq C_1\left( C_2\a_s+C_3  r^{\min(2,p)} \a_s^{\min(2,p)}\right).
$$
There exists   $r_0>0$ large    and $s_0\in(0,1)$ (possibly depending on $r_0$) such that for any
$s\in (s_0,1)$ we have
\begin{eqnarray*}
 r_0\a_{s_0}&\geq& C_1\left( C_2\a_{s_0} +C_3  r_0^{\min(2,p)} \a_{s_0}^{\min(2,p)}\right)\\
&>& C_1\left( C_2\a_s +C_3  r_0^{\min(2,p)} \a_{s_0}^{\min(2,p)}\right),
\end{eqnarray*}
since $\a_s$ is small as $s\nearrow1$. It follows that
for every $s\in(s_0,1)$,
the mapping $$\o\mapsto  -(\Phi_{s}'(0))^{-1}\left\{\Phi_{s}(0)+Q_{s}(\o)\right\}$$ maps $\calB_{r_0,{s_0}}$ into itself.
Increasing $s_0$ if necessary,  this map is a contraction on $\calB_{r_0,{s_0}}$ by \eqref{eq:contrac}.
Hence by the Banach fixed point theorem, for every $s\in(s_0,1)$,
 there exists a unique function $\o^{s}\in \calB_{r_0,{s_0}}$ solving the fixed point equation \eqref{eq:fixed}.
\QED

The set of pseudo-minimizers is given by $\{U_1+\o_s\,:\, \Phi_s(\o_s)=0,\,s\in(s_0,1)\}$.
We now prove uniqueness, up to translations,  of the minimizers for $\nu_s$ when $s$ is close to $1$ by showing that minimizers belong to such a set.

\subsection{Completion of the proof of
Theorem \ref{th:Uniq}}

Let $u^1_s$ and $u^2_s$ be two minimizers for $\nu_s$. We know that
they are symmetric under rotation, so we may and do assume that they
are both symmetric with respect to the origin of $\R^N$.
Our aim is to show that $u^1_s= u^2_s$
provided $s$ is close to $1$ (no confusion
should arise between the superscripts $1$ and $2$ and some exponents
that shall occur in the course of the proof).

By Lemma \ref{lem:reg}, we know that
$ u^i_s= U_1+\o^i_s$ with $\|\o^i_s\|_s\to0$ as $s\nearrow 1$, for $i=1,2$ and $\o^i$
is symmetric with respect to the origin  for $i=1,2$. Then we have $\Phi_s(\o^i_s)=0$ for $s$ close to 1
and thus by uniqueness (Proposition \ref{prop:finite}) we conclude that $\o^1_s=\o^2_s$.

\subsection*{Acknowledgements}

E. V. is supported by the ERC Grant 277749
Elliptic Pde's and Symmetry of Interfaces and Layers for Odd Nonlinearities
(EPSILON).
M. M. F. is supported the by the Alexander von Humboldt foundation.

   \end{document}